\documentclass[11pt]{amsart}
\usepackage{latexsym,amssymb,amsmath}
\usepackage[all,ps,cmtip]{xy} 
\newtheorem{theo}{Theorem}
\newtheorem{coro}{Corollary}
\newtheorem{prop}{Proposition}
\newtheorem{lemm}{Lemma}

\begin{document}

\def\ot{\otimes}
\def\we{\wedge}
\def\wec{\wedge\cdots\wedge}
\def\op{\oplus}
\def\ra{\rightarrow}
\def\lra{\longrightarrow}
\def\fso{\mathfrak so}
\def\cO{\mathcal{O}}
\def\cS{\mathcal{S}}
\def\fsl{\mathfrak sl}
\def\PP{\mathbb P}\def\PP{\mathbb P}\def\ZZ{\mathbb Z}\def\CC{\mathbb C}
\def\RR{\mathbb R}\def\HH{\mathbb H}\def\OO{\mathbb O}
\def\smc{\cdots }
\title{On rectangular Kronecker coefficients}
\author[L. Manivel]{Laurent Manivel}
\address{Institut Fourier,  
Universit\'e de Grenoble I et CNRS,
BP 74, 38402 Saint-Martin d'H\`eres, France}
\email{{\tt Laurent.Manivel@ujf-grenoble.fr}}

\begin{abstract} 
We show that rectangular Kronecker coefficients stabilize when the lengths of the sides 
of the rectangle grow, and we give an explicit formula for the limit values in terms of 
invariants of $\fsl_n$.
\end{abstract}

\maketitle


\section{Introduction}
Kronecker coefficients are structure constants for tensor products of 
irreducible complex representations of symmetric groups. More precisely, 
the irreducible representations of $\cS_n$, the group of permutations of $n$ symbols, 
are  indexed by partitions of $n$, see \cite[I.7]{mcd}. Denote by $[\lambda]$ the representation 
associated to a partition $\lambda$. Then we can decompose 
$$[\lambda]\otimes [\mu] = \bigoplus_{\nu} k_{\lambda,\mu,\nu}[\nu],$$
and a major unsolved problem is find a combinatorial rule for computing 
the Kronecker coefficients $k_{\lambda,\mu,\nu}$. 

Motivated by certain aspects of algebraic complexity theory (see \cite{blmw}
for an overview), we focus on a more specific problem. Suppose that $\mu=\nu=(d^n)$, 
the {\it rectangle partition} with $n$ parts equal to $d$. Can we compute
the corresponding Kronecker coefficients $k_{\lambda,(d^n),(d^n)}$, which we call 
rectangular Kronecker coefficients? Our main result is the following:

\begin{theo}
Let $\lambda=(dn-|\rho|,\rho)$ for some partition $\rho$. 
Consider the Kronecker coefficient $k_\rho(d,n)=k_{\lambda,(d^n),(d^n)}$. 
\begin{enumerate}
\item $k_\rho(d,n)$ is a symmetric, non decreasing function of $n$ and $d$;
\item for $d\ge |\rho|$, the Kronecker coefficient $k_\rho(d,n)$ is equal to 
$$k_\rho(n) := \dim S_\rho (\fsl_n)^{GL_n}.$$
\end{enumerate}
\end{theo}

A few words of explanation are in order: the general complex linear group $GL_n$ acts by conjugation 
on the Lie algebra $\fsl_n$ of traceless matrices. For any partition $\rho$ there is an induced 
action of $GL_n$ on the Schur power $S_\rho (\fsl_n)$, and $S_\rho (\fsl_n)^{GL_n}$ denotes the 
subspace of invariants of this action. (Note that we could replace $GL_n$ by $SL_n$.)  
Also, we have denoted by $|\rho|$ the sum of the parts of the partition $\rho$. In the sequel
we will denote by $\ell(\rho)$ its length, defined as the number of non zero parts. 

\section{Proof of the theorem}
We first use Schur-Weyl duality to translate the computations of Kronecker coefficients to the 
setting of $GL_m$-modules and their Schur powers, or equivalently, to a problem involving symmetric 
functions. By Schur-Weyl duality, we mean the following statement (see, e.g., \cite{ho}). 
For a complex vector space $V$, 
the $\cS_n\times GL(V)$-module $V^{\otimes n}$ decomposes as 
$$V^{\otimes n} = \bigoplus_{\substack{|\lambda|=n, \\ \ell(\lambda)\le \dim V}}[\lambda]\otimes S_\lambda V.$$
A straightforward consequence is that, for two vector spaces $U,V$ and three partitions 
$\lambda,\mu,\nu$ of the same integer $n$, the multiplicity of $S_\mu U\otimes S_\nu V$
inside $S_\lambda(U\otimes V)$ is equal to the Kronecker coefficient 
$k_{\lambda,\mu,\nu}$ (at least for $\ell(\mu)\le\dim V, \ell(\nu)\le\dim W$ and 
$\ell(\lambda)\le\dim U \dim V$). 

Suppose that $V=W\oplus L$ for some vector spaces $W,L$, with $\dim L=1$. We compute 
$S_\lambda(U\otimes V)$ in two different ways. First, we can use the previous remark and write
\begin{equation}\label{e1}
S_\lambda(U\otimes V) =  \bigoplus_{\mu,\nu}k_{\lambda,\mu,\nu} S_\mu U\otimes S_\nu (W\oplus L) 
 =  \bigoplus_{\substack{\mu,\nu, \theta \\ \nu\mapsto\theta}}k_{\lambda,\mu,\nu} S_\mu U\otimes S_\theta W
\otimes L^{|\nu|-|\theta|}.
\end{equation}
Here we have used a version of the {\it Pieri formula}, 
$$S_\nu (W\oplus L)  =  \bigoplus_{\theta \mid \nu\mapsto\theta}S_\theta W\otimes L^{|\nu|-|\theta|},$$
where the notation $\nu\mapsto\theta$ means that $\nu_i\ge\theta_i\ge \nu_{i+1}$ for all $i$. 
This formula is a special case of a more general formula which is also a consequence of 
Schur-Weyl duality: applying a Schur functor to a direct sum, one gets
$$S_\lambda(E\oplus F)=\bigoplus_{\substack{\alpha,\beta \\ |\alpha|+|\beta|=|\lambda|}}c_\lambda^{\alpha,\beta}
S_\alpha E\otimes S_\beta F,$$
where the {\it Littlewood-Richardson coefficient} $c_\lambda^{\alpha,\beta}$ can be defined as the 
multiplicity of $S_\lambda G$ inside the tensor product  $S_\alpha G\otimes S_\beta G$, for $G$ of 
large enough dimension. We can apply this formula to derive 
$$S_\lambda(U\otimes V)=S_\lambda(U\otimes W\oplus U\otimes L)=
\bigoplus_{\substack{\alpha,\beta \\ |\alpha|+|\beta|=|\lambda|}}c_\lambda^{\alpha,\beta}
S_\alpha (U\otimes W)\otimes S_\beta U\otimes L^{|\beta|}$$
and then, decomposing $S_\alpha (U\otimes W)$ with the help of Kronecker coefficients, 
\begin{equation}\label{e2}
S_\lambda(U\otimes V)
=\bigoplus_{\substack{\alpha,\beta,\rho,\theta \\ |\alpha|+|\beta|=|\lambda|}}c_\lambda^{\alpha,\beta}
k_{\alpha,\rho,\theta}S_\rho U\otimes S_\beta U\otimes S_\theta W\otimes L^{|\beta|}.
\end{equation}
Now we specialize these identities to the case where $\lambda=(d^n)$. Equating the coefficients, in (\ref{e1})
and (\ref{e2}), of $S_{(d^n)}U\otimes S_\theta W\otimes L^{nd-|\theta|}$, we deduce the following equality:
\begin{equation}\label{e4}
\sum_{\nu\mapsto\theta}k_{(d^n),(d^n),\nu}=\sum_{\alpha,\beta,\rho}
c_{(d^n)}^{\alpha,\beta}k_{\alpha,\rho,\theta}c_{(d^n)}^{\rho,\beta}.
\end{equation}
The Littlewood-Richardson coefficients in this formula are easy to understand: 
we have the following folklore lemma, which is an easy application of the Littlewood-Richardson
rule (and also a version of Poincar\'e duality for Grassmannians). 

\begin{lemm} 
If the Littlewood-Richardson $c_{(d^n)}^{\alpha,\beta}$ is non zero, then $\alpha$ and $\beta$ are
complementary partitions in the rectangle $d\times n$; that is, 
$$\alpha_i+\beta_{n-i+1}=d \qquad \forall i.$$
(In particular $\alpha_1, \beta_1\le d$ and $\alpha_{n+1}=\beta_{n+1}=0$.)
When this condition is satisfied, $c_{(d^n)}^{\alpha,\beta}=1.$
\end{lemm}

This implies that in formula (\ref{e4}), can only contribute to the right-hand side the 
partitions $\alpha,\beta,\rho$ contained in the rectangle $d\times n$, and such that $\alpha=\rho$ is 
complementary to $\beta$. We can thus write:
\begin{equation}\label{e3}
\sum_{\nu\mapsto\theta}k_{(d^n),(d^n),\nu}=\sum_{\alpha\subset d\times n}k_{\alpha,\alpha,\theta}.
\end{equation}
Observe that these two sums are of a quite different nature: the first one is taken over partitions 
$\nu$ of size $dn$; but the second one is taken over partitions $\alpha$ of the same size as $\theta$; 
and this size can be arbitrary. In particular, if we suppose that $|\theta|\le\min(d,n)$, then 
the condition that $\alpha\subset d\times n$ is empty. If we simply ask that $|\theta|\le d$, it amounts 
to the condition that the length of $\alpha$ does not exceed $n$, and we can rewrite the previous 
identity as
$$\sum_{\nu\mapsto\theta}k_{(d^n),(d^n),\nu}=\sum_{\ell(\alpha)\le n}k_{\alpha,\alpha,\theta}.$$
In particular the right hand side does not depend on $d\ge |\theta|$.

Now, write $\nu=(dn-|\rho|,\rho)$ for some partition $\rho$, and observe that the condition that 
$\nu\mapsto\theta$ means that $dn-|\rho|\ge\theta_1\ge\rho_1\ge\theta_2\ge\rho_2\cdots $. Otherwise
said, we ask that $dn-|\rho|\ge\theta_1$ and $\theta\mapsto\rho$.
In particular, the latter condition implies that, 
$|\rho|\le |\theta|$, and therefore the condition that $dn-|\rho|\ge\theta_1$ 
follows automatically as 
soon as $d\ge 2$. We can thus rewrite our identity as 
\begin{equation}\label{e5}
\sum_{\theta\mapsto\rho}k_{(d^n),(d^n),(dn-|\rho|,\rho)}=\sum_{\ell(\alpha)\le d}k_{\alpha,\alpha,\theta}.
\end{equation}
In order to interprete the right-hand side, we apply the formula for Schur powers of tensor products
to two vector spaces in duality. That is, we write
$$S_\theta(V\otimes V^\vee)=\bigoplus_{\substack{\alpha,\beta, \\ \ell(\alpha),\ell(\beta)\le\dim V}}
k_{\alpha,\beta,\theta}S_\alpha V\otimes S_\beta V^\vee.$$
By Schur's lemma, the space of $GL(V)$-invariants inside $S_\alpha V\otimes S_\beta V^\vee=Hom(S_\beta V,S_\alpha V)$  
is non zero if and only if $\alpha=\beta$, in which case it is one-dimensional. Letting $V={\mathbb C}^n$, we 
deduce that 
$$\dim S_\theta (End_n)^{GL_n} = \sum_{\ell(\alpha)\le n}k_{\alpha,\alpha,\theta}.$$
But $End_n=\fsl_n\oplus {\mathbb C}$, hence by the Pieri formula $S_\theta (End_n)=\oplus_{\theta\mapsto\rho}
S_\rho(\fsl_n)$. We can thus rewrite the identity (\ref{e5}) between Kronecker coefficients as 
$$\sum_{\theta\mapsto\rho}k_{(d^n),(d^n),(dn-|\rho|,\rho)}=\sum_{\theta\mapsto\rho}\dim S_\rho(\fsl_n)^{GL_n},$$
for any $\theta$ such that $|\theta|\le d$. But this implies that for any $\rho$ such that $|\rho|\le d$, 
$$k_{(d^n),(d^n),(dn-|\rho|,\rho)}=\dim S_\rho(\fsl_n)^{GL_n}.$$
Indeed, $\theta\mapsto\rho$ implies $\theta=\rho$ or $|\rho|<|\theta|$, so by induction we are immediately done. 
This proves the second assertion of the Theorem. 

\medskip
The proof of the first assertion is straightforward. First, the sign representation $[1^m]$ of $\cS_m$ has the 
property that its tensor product with an irreducible representation $[\lambda]$ is the irreducible representation 
$[\lambda^\vee]$, where $\lambda^\vee$ denotes the transpose partition to $\lambda$. Since $[1^m]\otimes [1^m]$
is just the trivial representation, this implies that $[\lambda]\otimes [\lambda]=[\lambda^\vee]\otimes [\lambda^\vee]$.
Applying this remark to a rectangular partition $\lambda=(d^n)$, we get $\lambda^\vee=(n^d)$, hence 
$$k_\rho(d,n)=k_\rho(n,d)\qquad\forall\rho.$$
Second, this is a non decreasing function of $d$ because of the semi-group property of Kronecker coefficients:
if $k_{\lambda,\mu,\nu}$ and $k_{\lambda',\mu',\nu'}$ are non zero, then  $k_{\lambda+\lambda',\mu+\mu',\nu+\nu'}$
is also non-zero, and more precisely, 
$$k_{\lambda+\lambda',\mu+\mu',\nu+\nu'}\ge \max(k_{\lambda,\mu,\nu}, k_{\lambda',\mu',\nu'}).$$
Indeed, let $E,F,G$ be three vector spaces and $U$ be the unipotent subgroup of $GL(E)\times
GL(F)\times GL(G)$ consisting of triples of strictly triangular matrices, with respect to some choice of
basis. Let $T$ be the maximal torus consisting of triples of diagonal matrices in the same basis. 
Then it is another consequence of Schur Weyl duality that the invariant algebra $A=Sym(E\otimes F\otimes G)^U$,
which is a $T$-modules, has weight spaces $A_{\lambda,\mu,\nu}$ of dimension $k_{\lambda,\mu,\nu}$. 
Then the semi-group property follows from the obvious fact that $A$ is an algebra without zero-divisors.  

Since $k_{1^n,1^n,n}$ is equal to one, we are done.  \qed

\medskip\noindent {\it Remark 1}. In the second assertion of the Theorem, we can improve the bound on $d$ as follows: 
$$k_\rho(d,n) = \dim S_\rho (\fsl_n)^{GL_n}$$
as soon as $2d\ge |\rho|+\rho_1$. Indeed,  for a Kronecker coefficient $k_{\alpha,\beta,\theta}$ to be non
zero, we need that $m-\theta_1\le (m-\alpha_1)+(m-\beta_1)$, where $m=|\alpha|=|\beta|=|\theta| $
\cite[Theorem 2.9.22]{JK}. Therefore $k_{\alpha,\alpha,\theta}\ne 0$ implies that $2\lambda_1\le |\theta|
+\theta_1$, so that in (\ref{e3}) we can replace the condition that $\alpha\subset d\times n$ by $\ell(\alpha)\le n$
as soon as $2d\ge |\theta|+\theta_1$.

\medskip\noindent {\it Remark 2}. The case where $n=2$ has been considered in \cite{man}, where we proved by 
completely different methods that $k_\rho(d,2)$ does not depend on $d$. More precisely, 
this Kronecker coefficient is equal to one
if and only if $\rho$ is even of length one or odd of length three (by an odd partition we mean a partition 
all of whose parts are odd), and zero otherwise. This result was also obtained independently in \cite{gwxz}. 
From the Theorem we deduce the corresponding statement 
for the dimension of $S_\rho(\fsl_2)^{GL_2}$.

\section{Complements}

\subsection{Invariants of matrices}
We have met  the identity 
$$\dim S_\theta (End_n)^{GL_n}=\sum_{\ell(\alpha)\le n}k_{\alpha,\alpha,\theta}.$$
For a given $\theta$, it implies that the dimension of $S_\theta (End_n)^{GL_n}$ is a non decreasing 
function of $n$, which is equal for $n\ge |\theta|=m$ to the multiplicity of $[\theta]$ inside
$$C_m = \bigoplus_{|\alpha|=m}[\alpha]\otimes [\alpha].$$
Observe that $C_m$ is nothing else that the conjugating representation of $\cS_m$ over $\mathbb{C}[\cS_m]$. 
Indeed, the regular representation  of $\cS_m$ over $\mathbb{C}[\cS_m]$, as a $\cS_m\times\cS_m$-module, 
is given by the same formula (but with exterior tensor products), and the conjugating representation is just the 
regular representation restricted to the diagonal embedding of $\cS_m$ inside  $\cS_m\times\cS_m$. 
We denote $C_m=\mathrm{diag}\; \mathbb{C}[\cS_m]$. 

A more transparent way to understand the previous identity is to observe that it is closely related 
to the fundamental theorems for invariants of matrices. Indeed (see, e.g., \cite[Chapter 11, \S 8]{procesi}, 
or \cite{fo}), 
the first fundamental theorem 
asserts that the $\cS_m$-equivariant map 
$$f_{m,n}:\mathrm{diag}\; \mathbb{C}[\cS_m]\longrightarrow (End_n^{\otimes m})^{GL_n}$$
mapping a permutation $\sigma$ to $f_{m,n}(\sigma)\in End_n^{\otimes m}=End(\mathbb{C}^{\otimes m})$
defined by 
$$f_{m,n}(\sigma)(X_1\otimes\cdots\otimes X_m)=X_{\sigma(1)}\otimes\cdots\otimes X_{\sigma(m)},$$
is surjective. Moreover, the second fundamental theorem asserts that the  kernel of $f_{m,n}$ is the sum of the 
submodules $[\alpha]\otimes [\alpha]$ of $\mathrm{diag}\; \mathbb{C}[\cS_m]$, for $\ell(\alpha)>n$. 
In particular, $f_{m,n}$
is an isomorphism for $n\ge m$.

\subsection{Derangements}
How can we pass from $End_n^{\otimes m}$ to $\fsl_n^{\otimes m}$? A first observation is that, since 
$End_n=\mathbb{C}I_n\oplus\fsl_n$, and the identity $I_n$ is certainly $GL_n$-invariant, we have
$$ (End_n^{\otimes m})^{GL_n}\simeq \bigoplus_{p=0}^m\binom{m}{p}(\fsl_n^{\otimes p})^{GL_n}.$$
Since the left hand side, for $n\ge m$, has dimension $m!$, we easily deduce that 
$$\dim (\fsl_n^{\otimes m})^{GL_n}=\sum_{p=0}^m(-1)^p \binom{m}{p}(m-p)!=D_m,$$
the number of derangements, that is, fixed-point free permutations in $\cS_m$.
This was already observed in  \cite{bd}, but only under the weaker condition   
that $n\ge 2m$. 

We can be more precise by making the following observation. We have
$$ (\fsl_n^{\otimes m})^{GL_n} = (End_n^{\otimes m})^{GL_n}\cap \fsl_n^{\otimes m},$$
and the subspace $\fsl_n^{\otimes m}$ of $End_n^{\otimes m}$ can be characterized as the intersection of
the kernels of the  contractions
$$c_k : End_n^{\otimes m} \ra End_n^{\otimes m-1}$$ 
defined, for $1\le k\le m$, by applying the trace morphism to the $k$-th factor of $End_n^{\otimes m}$.

Consider an element $f_{n,m}(\sigma)$ for some $\sigma\in\cS_m$. Once we have chose an basis $e_1,\ldots, e_n$,
we can write it explicitely as 
$$f_{n,m}(\sigma)=\sum_{i_1,\ldots ,i_m=1}^n e_{i_1}^\vee\otimes e_{i_{\sigma(1)}}\otimes
\cdots\otimes e_{i_n}^\vee\otimes e_{i_{\sigma(n)}}.$$
An easy computation then shows that a linear combination $a=\sum_\sigma a_\sigma \sigma\in \mathbb{C}[\cS_m]$ 
belongs to the kernel of $c_m$, the last contraction, if and only if 
$$na_{(\tau, m)}+\sum_{k=1}^{m-1}a_{s_{k,m} (\tau, m)}=0$$
for all $\tau\in\cS_{m-1}$. Here we have denoted by $(\tau, m)\in\cS_m$ the permutation deduced from 
$\tau$ by adding the fixed point $m$, and by $s_{k,m}$ the transposition exchanging $k$ and $m$.
In this identity, note that $(\tau, m)$ has one more fixed point than $\tau$, while $s_{k,m} (\tau, m)$
has either the same number of fixed points as $\tau$, or one less. This implies, inductively, 
that if $a=\sum_\sigma a_\sigma \sigma$ is annihilated by the contractions $c_1,\ldots ,c_m$, then 
all the coefficients of $a$ can be expressed in terms of the coefficients  $a_\sigma$, for $\sigma$
describing the set of fixed point free permutation -- or derangements. Since we know that the dimension 
of the space of such $a$'s is precisely the set of derangements, the remaining coefficients are 
independent. Let $\mathrm{diag}\; \mathbb{C}[\cS_m]_{fpf}\subset \mathrm{diag}\; \mathbb{C}[\cS_m]$
be the submodule defined by the derangements in $\cS_m$. We deduce from the preceding discussion:

\begin{prop}
For $n\ge m$, there is an isomorphism of $\cS_m$-modules
$$\mathrm{diag}\; \mathbb{C}[\cS_m]_{fpf}\simeq (\fsl_n^{\otimes m})^{GL_n}.$$
\end{prop}

This has the following consequence for rectangular Kronecker
coefficients. 

\begin{coro}
The limit value $k_\rho$ of $k_{\rho}(d,n)$ for $d,n\ge |\rho|$, is equal to the 
multiplicity of $[\rho]$ inside $\mathrm{diag}\; \mathbb{C}[\cS_m]_{fpf}$,
the fixed point free part of the conjugating representation.
\end{coro}
 
Although $D_m$ grows very fast with $m$, the coefficients $k_{\rho}$ seem to 
remain much smaller. For $m\le 6$ their values are given by the following  table:
$$\begin{array}{lllllll}
 m=2 &k_2=1,\; &k_{11}=0,\; & & & &D_2=1,\\
 m=3 &k_3=1,\; &k_{21}=0,\; &k_{111}=0,\; & & &D_3=2, \\
 m=4 &k_4=2,\; &k_{31}=0,\; &k_{22}=2,\; &k_{211}=1,\; &k_{1^4}=0,\; &D_4=9, \\
 m=5 &k_5=2,\; &k_{41}=1,\; &k_{32}=2,\; &k_{311}=3,\;&k_{221}=1,\;  & \\
  &k_{21^3}=1,\; &k_{1^5}=1, & & & &D_5=44, \\
 m=6,&k_6=3,\; &k_{51}=1,\; &k_{42}=6,\; &k_{33}=1,\;&k_{411}=4,\;  & \\
 &k_{321}=4, &k_{2^3}=5,\; &k_{31^3}=4,\; &k_{2211}=2,\; &k_{21^4}=2, \; \\
  & k_{1^6}=0, & & & &  &D_6=265.
\end{array}$$

\subsection{Symmetric and skew-symmetric invariants}
Certain submodules of the tensor algebra of $\fsl_n$ have a classical interpretation. 
Indeed, it is
well-known that $GL_n$-invariants in the symmetric algebra $Sym(\fsl_n)$ generate a polynomial
algebra with generators in degree $2,3,\ldots ,n$. On the other hand, $GL_n$-invariants in the 
exterior algebra $\Lambda (\fsl_n)$ can be interpreted as invariant differential forms on the 
Lie group $SL_n$, and it follows from classical results by Hopf and Samelson that 
$\Lambda (\fsl_n)^{GL_n}$ is an exterior algebra with generators of odd degrees $3, 5, \ldots , 2n-1$. 
Rephrasing these facts, and taking Remark 1 into account, we get:

\begin{coro}

\begin{enumerate}
\item
For $d\ge m$, the Kronecker coefficient $k_{(m)}(d,n)$ is equal to the number
of partitions of $m$ into integers between $2$ and $n$. 
\item
For $d\ge \frac{m+1}{2}$, the Kronecker coefficient $k_{(1^m)}(d,n)$ 
is equal to the number of partitions of $m$ into distinct odd integers between $3$ and $2n-1$. 
\end{enumerate}
\end{coro}

Note that partitions of $m$ into distinct odd integers are in correspondence with partitions 
$\lambda$ of $m$ equal to their own transpose. This is equivalent to the  condition that 
$[\lambda]\otimes [\lambda]$ contains the sign representation, so that the number of 
partitions of $m$ into distinct odd integers can be expressed as 
$$p_{odd}^0(m)=\sum_{|\lambda|=m}k_{\lambda,\lambda,(1^m)}.$$

\subsection{A birational map}
It may be possible to give a more transparent and more geometric proof of our results 
starting from the following observation. Let $E,F,G$ be three vector spaces of dimensions
$n,n,g$, with $g\le n^2$. Let $\gamma_1,\ldots ,\gamma_g$ be a basis of $G$, and denote by $U_G$ the 
unipotent subgroup of $GL(G)$ consisting of endomorphisms whose matrix in this basis
is strictly upper-triangular.  Consider a generic element $T\in E^\vee\otimes F\otimes G^\vee
=Hom(G,Hom(E,F))$. The image $T(\gamma_1)$ of $\gamma_1$ is an invertible morphism, which 
identifies $E$ with $F$. Under this identification, $T(\gamma_1)$ becomes the identity map
and we can deduce a tensor $\bar T\in Hom(G/\langle  \gamma_1\rangle), End(E)/\langle Id_E\rangle)$,
where $End(E)/\langle Id_E\rangle$ can be identified with $\fsl(E)$. Since $T$ is generic and 
$g\le n^2$, $\bar T$ is injective, and we deduce a flag $W_1\subset \cdots\subset W_{g-1}\subset
\fsl(E)$, where $W_k$ has dimension $k$. Denoting by $\mathcal{F}_{g-1}(\fsl(E))$ the variety 
parametrizing such flags, this construction defines a rational map
$$\pi: \PP(E^\vee\otimes F\otimes G^\vee)//SL(E)\times SL(F)\times U_G \longrightarrow
\mathcal{F}_{g-1}(\fsl(E))//SL(E)$$
which is easily seen to be birational. 
The left-hand side is $Proj(A)$, where 
$$A=\bigoplus_{\delta\ge 0}S^\delta(E^\vee\otimes F\otimes G^\vee)^{SL(E)\times SL(F)\times U_G}
=\bigoplus_{d\ge 0}k_{(d^n),(d^n),\lambda}(S_\lambda G^\vee)^{U_G}$$
is endowed with a natural action of the torus $T_G$ of diagonal matrices in $GL(G)$, such
that the one-dimensional vector space $(S_\lambda G^\vee)^{U_G}$ is acted on by the character
defined by $\lambda$. In consequence, $A$ is graded over the character group of $T_G$, and the dimension
of $A_\lambda$ is the rectangular Kronecker coefficient $k_{(d^n),(d^n),\lambda}$.

On the right-hand side, consider a line bundle $\mathcal{L}_\alpha$, whose fiber over a point 
of $\mathcal{F}_{g-1}(\fsl(E))$ defined by a flag $W_1\subset \cdots\subset W_{g-1}\subset
\fsl(E)$ is $W_1^{\alpha_1}\otimes (W_2/W_1)^{\alpha_2}\otimes\cdots\otimes (W_{g-1}/W_{g-2})^{\alpha_{g-1}}$. If
$\alpha$ is a partition, then by the Borel-Weil theorem 
$$H^0(\mathcal{F}_{g-1}(\fsl(E)),\mathcal{L}_\alpha^{-1})\simeq S_\alpha \fsl(E).$$
Moreover, any $SL(E)$-invariant section of  $\mathcal{L}_\alpha^{-1}$ defines, through $\pi$, a rational
function on $\PP(E^\vee\otimes F\otimes G^\vee)//SL(E)\times SL(F)\times U_G$. Multiplying this function 
by a $d$-th power of $\gamma_1^n\in S^nG\simeq \Lambda^nE\otimes \Lambda^nF^\vee\otimes S^nG\subset
S^n(E\otimes F^\vee\otimes G)$, this rational function will become regular for $d$ large enough 
(this is because in the construction of $\pi$ we have inverted $T(\gamma_1)$, and we are now 
multiplying by some power of its determinant). 
The weight of the resulting section with respect to $T_G$ is easily computed to be $(nd-|\alpha|,\alpha)$,
and we finally get a map 
$$S_\alpha \fsl(E)^{SL(E)}\longrightarrow A_{(nd-|\alpha|,\alpha)}.$$
Our theorem can thus be translated into the claim that, for $d$ large enough, this map
is an isomorphism. 


\subsection{A question}. 
We have seen that the semi-group property of Kronecker coefficients, and the obvious fact 
that $k_{(1^n),(1^n),(n)}=1$, is enough to ensure that the Kronecker coefficient we denoted
$k_{d,n}(\rho)$ is a non decreasing function of $n$ and $d$. More generally, for any integer
$\delta$ and any partition $\lambda$ such that $k_{(\delta^n),(\delta^n),\lambda}\ne 0$, 
we can conclude that for any integer $d$ and any partition $\mu$, the Kronecker coefficient 
$k_{(d+m\delta)^n,(d+m\delta)^n,\mu+m\lambda}$ is a non decreasing function of $m$. How
does it grow with $m$? Can one give an explicit equivalent when $m\ra +\infty$? When is there 
a finite limit? 

Of course this problem can also be considered for general Kronecker coefficients. 
The limit value of $k_{m+\alpha,m+\beta,m+\gamma}$, for $m\ra +\infty$, has been 
computed in \cite[3.4]{brion}.

\end{document}